\documentclass{article}
\usepackage{amsfonts}
\usepackage{amsmath}
\usepackage{amssymb}
\usepackage{amscd}

\newtheorem{theorem}{Theorem}

\newtheorem{definition}[theorem]{Definition}

\newtheorem{lemma}[theorem]{Lemma}

\newtheorem{proposition}[theorem]{Proposition}

\newcommand\R{\mathbb{R}}

\newcommand\K{\mathbb{K}}
\newcommand\C{\mathbb{C}}

\begin{document}

\title{n-Lie algebras }
\author{Michel GOZE - Nicolas GOZE - Elisabeth REMM \\
%EndAName
elisabeth.remm@uha.fr\\
Adress : LMIA.\ UHA\\
4 rue des Fr\`{e}res Lumi\`{e}re\\
F. 68093 Mulhouse Cedex}
\maketitle

\begin{abstract}
The notion of $n$-ary algebras, that is vector spaces with a multiplication concerning $n$-arguments, $n \geq 3$, became fundamental since the works of Nambu. Here we first present general notions concerning $n$-ary algebras and associative $n$-ary algebras. Then we will be interested in the notion of $n$-Lie algebras, initiated by Filippov, and which is attached to the Nambu algebras.
We study the particular case of nilpotent or filiform $n$-Lie algebras to obtain a beginning of classification.  This notion of $n$-Lie algebra admits a natural generalization in   Strong Homotopy $n$-Lie algebras in which the Maurer Cartan calculus is well adapted.
\end{abstract}

\medskip

\noindent {\small This work has been presented during the
 $11^{eme}$ {\it  Rencontre Nationale de G\'eom\'etrie Diff\'erentielle et Applications RNGDA'11,   Facult\'e des Sciences Ben-Msik, Casablanca.}
This meeting was dedicated to Professor  Younes BENSAID, Member of the Acad\'emie Fran\c{c}aise de Chirurgie.}

\medskip

\section{\protect\bigskip $n$-ary algebras}

\subsection{Basic definitions}

Let $\mathbb{K}$ be a commutative field of characteristic zero and $V$  a $%
\mathbb{K}$-vector space. Let $n$ be in \ $\mathbb{N}$, $n\geq 2.$ A $n$-ary
algebra structure on $V$ is given by a linear map
\begin{equation*}
\mu :V^{\otimes n}\rightarrow V.
\end{equation*}%
We denote by $(V,\mu )$ such an algebra. Classical algebras (associative
algebras, Lie algebras, Leibniz algebras for example) are binary that is
given by a $2$-ary product. In this paper, we are interested in the study of $%
n$-ary algebras for $n\geq 3.$ A
subalgebra of the $n$-ary algebra $(V,\mu )$ is a vector subspace $W$ of $V$
such that the restriction of $\mu $ to $W^{\otimes n}$ satisfies $\mu
(W^{\otimes n})\subset W.$ In this case $(W,\mu )$ is also a $n$-ary
algebra.

\begin{definition}
Let \ $(V,\mu )$ be a $n$-ary algebra.\ An ideal of $(V,\mu )$ is a
subalgebra $(I,\mu )$ satisfying%
\begin{equation*}
\mu (V^{\otimes p}\otimes I\otimes V^{\otimes n-p-1})\subset I,
\end{equation*}%
for all $p=0,\cdots ,n-1$ and where $V^{\otimes 0}\otimes I=I\otimes V^{\otimes
0}=I.$
\end{definition}

\begin{definition}
Let \ $(V_{1},\mu _{1})$ and $(V_{2},\mu _{2})$ be  $n$-ary algebras. A
morphism of $n$-ary algebras is a linear map $\varphi :V_{1}\rightarrow V_{2}
$ satisfying%
\begin{equation*}
\mu _{2}\circ \varphi ^{\otimes n}=\varphi \circ \mu _{1}.
\end{equation*}
\end{definition}

In this case, the linear kernel $Ker\varphi $ of the morphism $\varphi $ is
an ideal of $(V_{1},\mu _{1}).\ $In fact, $\ $if $v\in Ker\varphi $, then%
\begin{equation*}
\varphi (\mu _{1}(v_{1}\otimes \cdots \otimes v\otimes \cdots \otimes
v_{n-1}))=\mu _{2}(\varphi (v_{1})\otimes \cdots \otimes \varphi (v)\otimes
\cdots \otimes \varphi (v_{n-1}))=0.
\end{equation*}

To simplify notations, we identify the linear map $\mu$ on $V^{\otimes n}$ with the corresponding $n$-linear map on $V^n$. Then we write $\mu(v_1\otimes \cdots \otimes v_n)$ as well as  $\mu(v_1 \cdot v_2\cdots v_n)$.

\subsection{Anticommutative $n$-ary algebras}

Let \ $(V,\mu )$ be a $n$-ary algebra. It is called anticommutative if $\mu
(v_{1}\otimes \cdots \otimes v_{n})=0$ whenever  $v_{i}=v_{j}$ for $i\neq j.$
Since $\mathbb{K}$ is of characteristic $0$, this is equivalent to%
\begin{equation*}
\mu (v_{\sigma (1)}\otimes \cdots \otimes v_{\sigma (n)})=(-1)^{\varepsilon
(\sigma )}\mu (v_{1}\otimes \cdots \otimes v_{n}),
\end{equation*}%
for any $\sigma $ in the symmetric group $\Sigma _{n}$ where $%
(-1)^{\varepsilon (\sigma )}$ is the signum of the permutation $\sigma .$ If
$\mu $ is an antisymmetric $n$-ary multiplication, we write
\begin{equation*}
\lbrack v_{1},\cdots, v_{n}]
\end{equation*}%
instead of $\mu (v_{1}\otimes \cdots \otimes v_{n}).$

\subsection{Symmetric and commutative $n$-ary algebras}

A $n$-ary algebra $(V,\mu )$ is called symmetric  if it satisfies%
\begin{equation*}
\mu (v_{\sigma (1)}\otimes \cdots \otimes v_{\sigma (n)})=\mu (v_{1}\otimes
\cdots \otimes v_{n}),
\end{equation*}%
for all $v_{1},\cdots ,v_{n}$ $\in $ $V$ and for all $\sigma \in \Sigma _{n}.
$ It is called commutative if
\begin{equation*}
\sum_{\sigma \in \Sigma _{n}}(-1)^{\varepsilon (\sigma )}\mu
(v_{\sigma (1)}\otimes \cdots \otimes v_{\sigma (n)})=0,
\end{equation*}%
for all $v_{1},\cdots ,v_{n}$ $\in $ $V.$ Of course, any symmetric $n$-ary
algebra is commutative.

\subsection{Derivations}
Let $(V,\mu)$ a $n$-algebra.
\begin{definition}
A derivation of the $n$-algebra $(V,\mu)$ is a linear map
$$D: V\rightarrow V$$
satisfying
$$D(\mu(v_1,\cdots,v_n))=\sum _{i=1}^n \mu(v_1,\cdots,D(v_i),\cdots,v_n),$$
for any $v_1,\cdots,v_n \in V$.
\end{definition}
All derivations of $(V,\mu)$ generate a subalgebra of Lie algebra $gl(V)$. It is called the algebra of derivations of $V$ and denoted by $Der(V)$.

\medskip

\noindent {\bf Remark.} For any $v_1,\cdots,v_{n-1}$ in $V$, let $ad(v_1,\cdots,v_{n-1})$ be the linear map given by
$$ad(v_1,\cdots,v_{n-1})(v)=\mu(v_1,\cdots,v_{n-1},v).$$
Then this linear map is a (inner) derivation if and only if the product $\mu$ satisfies
\begin{eqnarray}\label{nLeibniz}
\mu(v_1,\cdots,v_{n-1},\mu(u_1,\cdots,u_n))=\sum_{i=1}^n
\mu(u_1,\cdots,\mu(v_1,\cdots,v_{n-1},u_i),\cdots,u_n).
\end{eqnarray}
We will study such a product for $n$-Lie algebras. If $n=2$, this shows that the maps $ad(X)$ are derivations of $(V,\mu)$ if and only if the binary product satisfies
$$\mu(v_1,\mu(u_1,u_2))=\mu(\mu(v_1,u_1),u_2)+\mu(u_1,\mu(v_1,v_2))$$
and $(V,\mu)$ is a Leibniz algebra (\cite{Cu}). Thus, for any $n$, a $n$-algebra
$(V,\mu)$ satisfying Equation (\ref{nLeibniz}) is called  $n$-Leibniz algebra.

\subsection{Simple, nilpotent $n$-ary algebras}
\begin{definition}
A $n$-ary algebra $(V,\mu)$ is called simple if
\begin{itemize}
\item $\mu$ is not abelian (i.e $\mu(V\cdots V) \neq 0$).
\item Any ideal is isomorphic to $V$ or is equal to ${0}$.
\end{itemize}
\end{definition}
We define the derived series by
$$
\left\{
\begin{array}{l}
V^{(1)}=V,\\
V^{(k)}=\mu(V^{(k-1)},V^{(k-1)},V,\cdots,V),
\end{array}
\right.
$$
and the lower central series by
$$
\left\{
\begin{array}{l}
V^{1}=V,\\
V^{k}=\mu(V^{k-1},V,V,\cdots,V).
\end{array}
\right.
$$
\begin{definition}
The $n$-ary algebra $(V,\mu)$ is called
\begin{itemize}
\item \emph{Solvable} if there is an integer $k$ such that $V^{(k)}=0.$
\item \emph{Nilpotent} if there is an integer $k$ such that $V^k=0.$
\end{itemize}
\end{definition}
The definitions presented here are the definitions given in \cite{Kas}.

\section{Particular classes of $n$-ary algebras}

\subsection{Associative $n$-ary algebras}

The extension of the notion of associativity to $n$-ary product is not
evident. The most natural is perhaps the notion of totally associativity
(see definition below).  Unfortunately, this notion is not auto-dual in the
 operadic point of view.\ This means that it is necessary to
introduce a most general notion of associativity called the partially
associativity.

\begin{definition}
The $n$-ary algebra $(V,\mu )$ is

$\bullet $ partially associative if $\mu $ satisfies
\begin{equation}
\displaystyle{\sum_{p=0}^{n-1}(-1)^{p(n-1)}\mu \circ (I_{p}\otimes \mu
\otimes I_{n-p-1})=0,}
\end{equation}

$\bullet $ totally associative if $\mu $ satisfies
\begin{equation}
\displaystyle{\ \mu \circ (\mu \otimes I_{n-1})=\mu \circ (I_{p}\otimes \mu
\otimes I_{n-p-1})},
\end{equation}%
for any $p=0,\cdots ,n-1.$
\end{definition}

\subsection{$\protect\sigma -$Associative $n$-ary algebras}

We found for the first time this notion in the study of symmetric spaces with
the Jordan point of view in the work of Loos. Recently, we have
studied a natural generalization of the matrices product.\ This example
gives directly the following definition:

Let $\sigma $ be an element of the symmetric group $\Sigma _{n}$ of order $n$%
. It induces an endomorphism of $V^{\otimes n}$ considering
\begin{equation*}
\Phi _{\sigma }^{V}(e_{i_{1}}\otimes \cdots \otimes e_{i_{n}})=e_{i_{\sigma
^{-1}(1)}}\otimes \cdots \otimes e_{i_{\sigma ^{-1}(n)}}.
\end{equation*}

\begin{definition}
The $n$-ary algebra $(V,\mu)$ is

$\bullet$ $\sigma$-partially associative if $\mu$ satisfies
\begin{equation}  \label{sigmapar}
\displaystyle{\sum_{p=0}^{n-1} (-1)^{p(n-1)} (-1)^{np\varepsilon(\sigma)}
\mu \circ (I_p \otimes (\mu \circ \Phi_{\sigma^{np}}^V) \otimes I_{n-p-1})=0}
\end{equation}

$\bullet$ $\sigma$-totally associative if $\mu$ satisfies
\begin{equation}  \label{sigmatot}
\ \mu \circ (\mu \otimes I_{n-1}) = \mu \circ (I_p \otimes
(\mu \circ \Phi_{\sigma^{np}}^V) \otimes I_{n-p-1}),
\end{equation}
for any $p=0, \cdots, n-1.$
\end{definition}

In particular, if $\sigma =Id$, we find the notion of partially and totally
associativity again.

\medskip

\noindent \textbf{Examples.}

\medskip

$\bullet $ \textbf{The Gerstenhaber products.} Let $\mathcal{A}$ be a
(binary) associative algebra and $H^{\ast }(\mathcal{A},\mathcal{A})$ its
Hochschild cohomology. The space of $k$-cochains is $\mathcal{C}^{k}(%
\mathcal{A})=Hom_{\mathbb{K}}(\mathcal{A}^{\otimes k},\mathcal{A}).$
A graded pre-Lie algebra $\oplus _{k}\mathcal{C}^{k}(%
\mathcal{A})$ has been defined by Gerstenhaber considering the product
\begin{equation*}
\bullet _{n,m}:\mathcal{C}^{n}(\mathcal{A})\times \mathcal{C}^{m}(\mathcal{A}%
)\rightarrow \mathcal{C}^{n+m}(\mathcal{A})
\end{equation*}%
given by
$$
\begin{array}{l}
\medskip
(f\bullet _{n,m}g)(X_{1}\otimes \cdots \otimes
X_{n+m-1})=\\
 \ \ \  \sum_{i=1}^{m}(-1)^{(i-1)(m-1)}f(X_{1}\otimes \cdots \otimes
g(X_{i}\otimes \cdots X_{i+m-1})\otimes \cdots \otimes X_{n+m-1}).
\end{array}
$$
A $k$-cochain $\mu $ satisfying $\mu \bullet _{k,k}\mu =0$ provides $%
\mathcal{A}$ with a $k$-ary partially associative structure.

\medskip

$\bullet$ \textbf{Some $n$-ary products of hypercubic matrices}

These products are defined in \cite{NG.R}. We consider the space $%
T_{q}^{p}(E)$ of tensors of $q$ contravariant and $p$ covariant type on a
vector space $E$. If $\left\{ e_{1},\cdots ,e_{n}\right\} $ is a basis of $E$
and $\varphi \in T_{q}^{p}(E)$, then
\begin{equation*}
\varphi (e_{i_{1}}\otimes \cdots \otimes e_{i_{p}})=\sum_{j_{1},\cdots
,j_{q}}C_{i_{1},\cdots ,i_{p}}^{j_{1},\cdots ,j_{q}}e_{j_{1}}\otimes \cdots
\otimes e_{j_{q}}.
\end{equation*}%
We denote by $\widetilde{\varphi }$ the tensor of $T_{p}^{q}(E)$ given by
\begin{equation*}
\widetilde{\varphi }(e_{j_{1}}\otimes \cdots \otimes
e_{j_{q}})=\sum_{i_{1},\cdots ,i_{p}}C_{i_{1},\cdots ,i_{p}}^{j_{1},\cdots
,j_{q}}e_{i_{1}}\otimes \cdots \otimes e_{i_{p}}.
\end{equation*}%
The $(2k+1)$-ary multiplication on $T_{q}^{p}(E)$ given by
\begin{equation*}
\mu (\varphi _{1},\cdots ,\varphi _{2k+1})=\varphi _{1}\circ \widetilde{%
\varphi _{2}}\circ \varphi _{3}\circ \cdots \circ \widetilde{\varphi _{2k}}%
\circ \varphi _{2k+1}
\end{equation*}%
is $s_{k}$-totally associative where $s_{k}\in \Sigma _{2k+1}$ is the
permutation
\begin{equation*}
s_{k}(1,\cdots ,2k+1)=(2k+1,2k,\cdots ,k,\cdots ,2,1).
\end{equation*}%
Identifying $\varphi $ to its structural constants $(C_{i_{1},\cdots
,i_{p}}^{j_{1},\cdots ,j_{q}})$ and to the hypercubic matrix $%
(C_{i_{1},\cdots ,i_{p},j_{1},\cdots ,j_{q}})$ we obtain a "natural"
extension of the classical product of matrices.

\section{$n$-Lie algebras}

Many notions of $n$-Lie algebras have been presented to generalize Lie algebras for $n$-ary algebras. The first one is probably due to Filippov (\cite{Fi}). These algebras have been studied from an algebraic point of view (classification, simplicity, nilpotency, representations) and because of their relations with the Nambu mechanic. The second one is the notion introduced with the strong homotopy algebra point of view. In this paper we are concerned by the two approaches.  To distinguish these different definitions we will call $n$-Lie algebras the first one and Lie $n$-algebras or sh-$n$-Lie algebras in the second one. In this section, we study  Filippov algebras.

\begin{definition}
An anticommutative $n$-ary algebra is a $n$-ary Lie algebra or simpler $n$%
-Lie algebra if the following Jacobi identity holds:%
\begin{eqnarray*}
\lbrack \lbrack u_{1},\cdots ,u_{n}],v_{1},\cdots ,v_{n-1}]
= \sum_{i=1}^n[u_{1},\cdots ,u_{i-1},[u_i,v_{1},\cdots
,v_{n-1}],u_{i+1},\cdots ,u_{n}],
\end{eqnarray*}%
for any $u_{1},\cdots ,u_{n},v_{1},\cdots ,v_{n-1}\in V.$
\end{definition}
This last condition is called Jacobi identity for $n$-Lie algebras.

\subsection{Fundamental examples}

1. This example was given by Fillipov. Let $A$ be a $n$-dimensional vector space on $\K$. Let $\{v_1,\cdots,v_{n+1}\}$ be a basis of $V$. The following product
$$[v_1,v_2,\cdots,\hat{v_i},\cdots,v_{n+1}]=(-1)^{n+1+i}v_i,$$
for $i=1,\cdots,n+1$
provides $A$ with a $n$-Lie algebra structure. We denote this algebra $A_{n+1}$.
\begin{theorem}
If $\K=\C$,  every simple $n$-Lie algebra is of dimension $n+1$ and it is isomorphic to
$A_{n+1}$.
\end{theorem}

\noindent 2. Let $A=\K[X_1,\cdots,X_n]$ be the associative algebras of $n$ indeterminates polynomials. We consider the product
$$[P_1,\cdots,P_n]=Jac(P_1,\cdots,P_n),$$
where $Jac$ denotes the Jacobian, that is the determinant of the Jacobian matrix of partial derivatives of $P_1,\cdots,P_n$. Provided with this product, $A$ is an infinite dimensional $n$-Lie algebra.

\medskip

\noindent 3. The Nambu brackets. It generalizes directly the previous example. Let $A={\mathcal{C}}^{\infty}(\R^3)$ be the algebra of differential functions on $\R^3$. This algebra is considered as classical observables on the three dimensional space $\R^3$ with coordinates $x,y,z$. We consider on $A$ the $3$-product
$$\{f_1,f_2,f_3\}=Jac(f_1,f_2,f_3).$$
This product is a $3$-Lie algebra product which generalizes the usual Poisson bracket from binary to ternary  operations.

\medskip

\noindent 4. (\cite{williams}). Let $A=\K[X_1,\cdots,X_n]$ be the associative algebra of $n$ indeterminates polynomials. Let $I_r$ be the linear subspace of $A$ linearly generated by the monomials of $A$ of degree greater than or equal to $r$. Clearly $I_r$ is a subspace of $I_3$ as soon as $r \geq 3.$ We define $J_r=I_3/I_r$ for $r>3$. For any $Q_1,\cdots,Q_n \in J_r$ we put
$$[Q_1,\cdots,Q_n]=Jac(Q_1,\cdots,Q_n).$$
This product is a $n$-Lie algebra bracket and $Q$ is a finite dimensional nilpotent $n$-Lie algebra.

\medskip

\noindent 5. Every $n$-Lie algebra of dimension $n$

\begin{itemize}
\item with $n$  odd, is abelian, that is,
\begin{equation*}
\lbrack e_{1},\cdots ,e_{n}]=0,
\end{equation*}

\item with $n$  even, is abelian or isomorphic to
\begin{equation*}
\lbrack e_{1},\cdots ,e_{n}]=e_{1},
\end{equation*}%
where $\left\{ e_{1},\cdots ,e_{n}\right\} $ is a basis of $V.$
\end{itemize}
In fact, we can write
\begin{equation*}
\lbrack e_{1},\cdots ,e_{n}]=\sum a_{i}e_{i},
\end{equation*}%
and if this product is not zero, we can consider that%
\begin{equation*}
\lbrack e_{1},\cdots ,e_{n}]=e_{1}.
\end{equation*}%
Thus
\begin{eqnarray*}
\lbrack \lbrack e_{1},\cdots ,e_{n}],e_{2},\cdots ,e_{n}]
&=&(-1)^{n-1}e_{1}.
\end{eqnarray*}
But the Jacobi identity implies
\begin{eqnarray*}
\lbrack \lbrack e_{1},\cdots ,e_{n}],e_{2},\cdots ,e_{n}]
&=&[[e_{1},e_{2},\cdots ,e_{n}],e_{2},\cdots
,e_{n}]+\cdots \\
&& +[e_{1},\cdots ,[e_{n},e_{2},\cdots ,e_{n}]]\\
&=&e_1.
\end{eqnarray*}
We obtain
$$(-1)^{n-1}e_{1}=e_1.$$

\subsection{Nilpotent $n$-Lie algebras}
In the first section we have defined nilpotency for general $n$-ary algebras. Since any $n$-Lie algebra is a $n$-Leibniz algebra, any adjoint operator
$ad(v_1,\cdots,v_{n-1})$ is a derivation.
\begin{theorem}(\cite{Kas})
For any finite dimensional nilpotent Lie algebras, the adjoint operators are nilpotent. Conversely, if the adjoint operators of the $n$-Lie algebra $V$ are nilpotent, then $V$ is nilpotent.
\end{theorem}
Assume that $V$ is a finite dimensional complex nilpotent $n$-Lie algebra. We will generalize the notion of characteristic sequence of Lie algebras to $n$-Lie algebras. We consider the set of generators of $V$ which is isomorphic to $V/V^2$.
\begin{lemma}
$$\dim V/V^2 \geq n.$$
\end{lemma}
Let us consider a free family $\{v_1,\cdots,v_{n-1}\}$ of $n-1$ vectors of $V-V^2$. The operator $ad(v_1,v_2,\cdots,v_{n-1})$ is a linear nilpotent operator of $V$ admitting $v_1,\cdots,v_{n-1}$ as eigenvectors associated to the eigenvalue $0$. We consider now the ordered sequence of the similitude invariants (the dimensions of Jordan blocks ) of this operator. It is of type $(c_1,\cdots,c_k,1,\cdots,1)$ with at least $n-1$ invariant equal to $1$, corresponding to the dimension of the eigenspace generated by the eigenvectors $v_i$. We assume that $c_1\geq \cdots\geq c_k\geq 0.$ We denote this sequence $c(v_1,\cdots,v_{n-1})$.
\begin{definition}
The characteristic sequence of the nilpotent $n$-Lie algebra is the sequence
$$c(V)=max\{c(v_1,\cdots,v_{n-1})\},$$
where $(v_1,\cdots,v_{n-1})$ are $n-1$ independent vectors of $V-V^2$, the order relation being the lexicographic order.
\end{definition}
Assume that $\dim V=p$. The possible extremal values of $c(V)$ are
 \begin{itemize}
 \item $(1,\cdots,1)$ and $V$ is an abelian $n$-Lie algebra,
 \item $(p-n+1,1,\cdots,1)$. This sequence  corresponds to a nilpotent operator $ad(v_1,v_2,\cdots,v_{n-1})$ with a biggest nilindex.
 \end{itemize}
\begin{definition}
A $p$-dimensional complex (or real) nilpotent $n$-Lie algebra is called filiform is its characteristic sequence is equal to $(p-n+1,\underbrace{1,\cdots,1}_{n-1})$.
\end{definition}
{\bf Examples}
\begin{itemize}
\item
We consider $n=3$ and $p=4$. The characteristic sequence is equal to $(2,1,1)$. Let $\{v_1,v_2,v_3,v_4\}$ be  a basis of $V$ such that the characteristic sequence of $ad(v_1,v_2)$ is $(2,1,1)$. If $\{v_3v_4\}$ is the Jordan basis of this operator then we have
$$[v_1,v_2,v_3]=v_4.$$
From the classification of \cite{Bai-Song}, we deduce that we have obtained the only filiform $3$-Lie algebra of dimension $4$.
\item We generalize easily this example. Let $V$ be the $p$-dimensional $3$-Lie algebra given by
    $$
    \left\{
    \begin{array}{l}
    [X_1,X_2,X_3]=X_4, \\
    \lbrack X_1,X_2,X_4 \rbrack=X_5, \\
    \cdots \\
     \lbrack X_1,X_2,X_{p-1}\rbrack =X_p.
    \end{array}
    \right.
    $$
    It is also a filiform $3$-Lie algebra. It is a model (\cite{GoM}) of the filiform $3$-Lie algebras of dimension $p$, that is any $p$-dimensional filiform $3$-Lie algebras can be contracted on this algebra.
    \item Every filiform $5$-dimensional $3$-Lie algebra is isomorphic to
    $$
    \left\{
    \begin{array}{l}
    [X_1,X_2,X_3]=X_4, \\
    \lbrack X_1,X_2,X_4\rbrack=X_5, \\
    \lbrack X_1,X_3,X_4\rbrack=aX_5, \\
    \lbrack X_2,X_3,X_4\rbrack=bX_5.\\
    \end{array}
    \right.
    $$
    \end{itemize}

\subsection{Graded filiform $n$-Lie algebras}
Let $f$ be a derivation of a complex filiform $n$-Lie algebra $V$. We assume that $f$ is diagonalizable. The decomposition of eigenspaces of $V$ gives a graduation of this $n$-Lie algebra. We consider the maximal abelian subalgebra of $Der(V)$ given by the diagonalizable derivations of $V$ which commute with $f$. We denote this algebra $T(f)$.
\begin{definition}
The rank of $V$ is the biggest dimension amongst the dimensions of $T(f)$ for any diagonalizable derivation $f$.
\end{definition}
\begin{proposition}
The rank of any filiform $n$-Lie algebra is equal to or smaller than $n$.
\end{proposition}
{\it Proof.} We consider the model given by
$$[X_1,X_2,\cdots,X_{n-1},X_i]=X_{i+1},$$
for $i=n+1,\cdots,p-1$, with $p=\dim V$. We can assume that $X_1,X_2,\cdots,X_n$ are eigenvectors. If we put
$$f(X_t)=\lambda _t X_t,$$
for $t=1,\cdots,n$, then other eigenvalues are
$$\lambda _i = \lambda _1+ \cdots +\lambda _{n-1}+\lambda_{i-1}$$
and this implies
$$\lambda _i = (n-i)(\lambda_1 + \cdots + \lambda_{n-1})+\lambda _n.$$
Thus $\lambda_1,\cdots, \lambda_n$ are the independent roots of this algebra which is then of rank $n$. Let $V_1$ be any filiform $n$-Lie algebra of dimension $p$. There exists $(X_1,\cdots,X_{n-1})$ independent vectors in $V_1-V_1^2$ such that the characteristic sequence of $V_1$ is given by the nilpotent operator $ad(X_1,\cdots,X_{n-1})$. We consider the corresponding Jordan basis of $V_1$. It satisfies
$$[X_1,X_2,\cdots,X_{n-1},X_i]=X_{i+1}$$
and other brackets are linear combinations of ${X_{n+1},\cdots,X_p}$. Let $f_t$ be the endomorphism given by $f_t(X_l)=X_l$ if $1\leq l\leq n$ and $f_t(X_l)=tX_l$ for $n+1\leq l \leq p.$ This endomorphism generates a contraction of $V_1$ in the model $V$. We deduce that the rank of $V_1$ is smaller than the rank of $V$.

\medskip

\begin{itemize}

\item Let us consider the filiform $3$-algebra
$$
    \left\{
    \begin{array}{l}
    [X_1,X_2,X_3]=X_4, \\
    \lbrack X_1,X_2,X_4\rbrack=X_5, \\
    \lbrack X_1,X_3,X_4\rbrack=aX_5, \\
    \lbrack X_2,X_3,X_4\rbrack=bX_5.\\
    \end{array}
    \right.
    $$
    Its rank is equal to $2$. In fact, in the basis $\{X_1,X_2,X_3-aX_2,X_4,X_5\}$, the algebra writes
    $$
    \left\{
    \begin{array}{l}
    [X_1,X_2,X_3]=X_4, \\
    \lbrack X_1,X_2,X_4\rbrack=X_5, \\
    \lbrack X_1,X_3,X_4\rbrack=0, \\
    \lbrack X_2,X_3,X_4\rbrack=bX_5.\\
    \end{array}
    \right.
    $$
    Any diagonalizable derivation which admits this basis as eigenvectors basis, satisfies
    $$f(X_i)=\lambda_i X_i$$
   with
   $$\lambda_3=\lambda_1, \lambda_4=2\lambda_1+\lambda_2,\lambda_5=3\lambda_1+2\lambda_2.$$
  Then the rank is $2$.

  \item For $n=2$, we have the following important result: any Lie algebra which admits a nonsingular derivation is nilpotent. This is false as soon as $n\geq 3$. We have the interesting example (\cite{williams}): consider the $n$-Lie algebra given by
      $$[X_1,X_2,\cdots,X_n]=X_2.$$
      This algebra admits a non singular derivation but it is not nilpotent.

      \item In a forthcoming paper we develop the classification of filiform
   $3$-Lie algebras whose rank is not $0$.

   \end{itemize}
\section{sh-$n$-Lie algebras or Lie $n$-algebras}

\subsection{Definition}
\begin{definition}
Let $\mu$ be a  $n$-ary skewsymmetric product on a vector space $A$. We say that
$(A,\mu)$ is a sh-$n$-Lie algebra (or a Lie $n$-algebra) if $\mu$ satisfies the (sh)-Jacobi's identity:
$$
\displaystyle \sum_{\sigma \in Sh(n,n-1)}(-1)^{\epsilon(\sigma)}\mu(\mu(x_{\sigma(1)},\cdots,x_{\sigma(n)}),
x_{\sigma(n+1)},
\cdots,x_{\sigma(2n-1)})=0,
$$
for any $x_i \in A$, where $Sh(n,n-1)$ is the subset of $\Sigma_{2n-1}$ defined by:
$$Sh(n,n-1)=\{\sigma \in \Sigma_{2n-1},\sigma(1)<\cdots<\sigma(n),\sigma(n+1)<\cdots<\sigma(2n-1)\}.$$
\end{definition}
Moreover, we assume that $\mu$ is of degree $n-2$.

\medskip

\noindent For example, if $n=3$, we have the following (sh)-Jacobi's identity, writing $(123)45$ in place of $\mu(\mu(x_1,x_2,x_3),x_4,x_5)$:
$$
\begin{array}{l}
(123)45-(124)35+(125)34+(134)25-(135)24+(145)23-(234)15+(235)14\\
-(245)13+(345)12=0.
\end{array}
$$

\subsection{$n$-Lie algebras and sh-$n$-Lie algebras}
\begin{proposition}
Any $n$-Lie algebra is a sh-$n$-Lie algebra.
\end{proposition}
{\it Proof. } The Jacobi condition for $n$-Lie algebras writes
$$\mu \circ (\mu \otimes I_{n-1}) \circ \Phi _v =0,$$
where $v \in \K[\Sigma _{2n-1}]$, the algebra group of the symmetric group  $\Sigma_{2n-1}$ on $2n-1$ elements, given by
$$
\begin{array}{ll}
v= &Id +\displaystyle \sum_{i=1}^n(-1)^i(i,n+1,\cdots,2n-1,1,2,\cdots,i-1,\widehat{i},i+1,\cdots,n),
\end{array}$$
where $(i,n+1,\cdots,2n-1,1,2,\cdots,i-1,\widehat{i},i+1,\cdots,n)$ is the permutation
$$
\left(
\begin{array}{lllllllllll}
1&2  &\cdots&n&n+1&n+2&\cdots&\cdots&\cdots&\cdots&2n-1\\
i&n+1&\cdots&2n-1&1&2&\cdots&i-1&i+1&\cdots&n
\end{array}
\right).
$$
Let
$$w=\sum_{\sigma \in \Sigma_{2n-1}}(-1)^{\epsilon(\sigma)}\sigma.$$
We have in $\K[\Sigma _{2n-1}]$, $w \circ v=\alpha(n)w$ with $\alpha(n)=1-n$ if $n$ is odd and $\alpha(n)=1+n$ if $n$ is even. This shows that the vector $w$ is in the invariant subspace of $\K[\Sigma _{2n-1}]$ generated by the vector $v$. This means that the identity
$$\mu \circ (\mu \otimes I_{n-1}) \circ \Phi _v =0$$
implies
$$\mu \circ (\mu \otimes I_{n-1}) \circ \Phi _w =0$$ which is equivalent to the Jacobi identity for sh-$n$-Lie algebras.
\begin{proposition}
A sh-$n$-Lie algebra is a $n$-Lie algebra if and only if any adjoint operator is a derivation.
\end{proposition}
{\it Proof. } We have seen that a $n$-Lie algebra is a $n$-Leibniz algebras and these last are characterized by the fact that any adjoint operator is a derivation.

\medskip

\noindent{\bf Remark. Colored Lie algebras, colored $n$-Lie algebras. } Let us consider a binary algebra with a skewsymmetric product satisfying a colored Jacobi identity:
$$\alpha[[X_i,X_j],X_k]+\beta[[X_j,X_k],X_l+\gamma[[X_k,X_i],X_j]=0,$$
for any $i<j<k$, the constants $\alpha,\beta,\gamma$ being in $\K$.
This identity is related to the vector $v=\alpha Id+\beta c+\gamma c^2$ of $\K[\Sigma_3]$. Let $w=Id-\tau_{12}-\tau_{13}-\tau_{23}+c+c^2$ the vector of $\K[\Sigma_3]$. Since $\K$ is of characteristic $0$, the Jacobi identity, is equivalent to
$$\mu \circ (\mu \otimes Id)\circ \Phi _w =0.$$
But in $\K[\Sigma_3]$ we have
$$w \circ v= (\alpha+\beta+\gamma)w.$$
Then, if $\alpha+\beta+\gamma \neq 0$, the colored Lie algebra satisfies the (non colored) Jacobi condition. This minimizes the interest of the notion of colored Lie algebras. It is the same for colored $n$-Lie algebras.

\subsection{$3$-Lie admissible algebras}

To simplify notations, we consider the case $n=3$. In this case the  product is of degree $1$. A $3$-ary algebra $(A,\cdot)$ is called $3$-Lie admissible if the antisymmetric product
$$[v_1,v_2,v_3]=\sum _{\sigma \in \Sigma _3}(-1)^{\varepsilon (\sigma)}v_{\sigma(1)}\cdot v_{\sigma(2)}\cdot v_{\sigma(3)}$$
is a sh-$3$-Lie product.
\begin{proposition}
A $n$-ary algebra $(A,\cdot)$ is  $3$-Lie admissible if and only if we have
$$
\begin{array}{l}
\displaystyle\sum _{\sigma \in \Sigma _5}(-1)^{\varepsilon (\sigma)}((v_{\sigma(1)}\cdot v_{\sigma(2)}\cdot v_{\sigma(3)})\cdot v_{\sigma(4)}\cdot v_{\sigma(5)}+
v_{\sigma(1)}\cdot (v_{\sigma(2)}\cdot v_{\sigma(3)}\cdot v_{\sigma(4)})\cdot v_{\sigma(5)}\\
+(v_{\sigma(1)}\cdot v_{\sigma(2)}\cdot (v_{\sigma(3)}\cdot v_{\sigma(4)}\cdot v_{\sigma(5)}))=0,
\end{array}$$
for any $v_1,v_2,v_3,v_4,v_5 \in A.$
\end{proposition}
{\bf Exemples.}
\begin{itemize}
\item Any $3$-ary partially associative algebra is $3$-Lie admissible.
\item In \cite{NG.R}, a notion of $\sigma$-associative algebra have been introduced in the space of tensors $(2,1)$ based on a vector space. In case of symmetric tensor, this product can be simplified. A symmetric tensor is defined by its structure constants $T_{ijk}$ which satisfy
$$ T_{ijk} = T_{jki} = T_{kij}. $$
The $3$-product $T \cdot U \cdot V$ whose structure constants are
$$(T \cdot U \cdot V)_{ijk}=\sum _l T_{lij}U_{lki}V_{ljk}$$
is $3$-Lie admissible. Moreover the associated sh-3-Lie algebra is a $3$-Lie algebra.
\end{itemize}

\subsection{Maurer-Cartan equations}

We assume in this section that any $n$-Lie algebras is of finite dimension. To simplify the presentation, we assume also that $n=3$. Let $V$ be a finite dimensional $3$-Lie algebras. Let $\{v_1,\cdots,v_p\}$ be a basis of $V$. The structure constants of $V$ related to this basis are given by
$$\{v_i,v_j,v_k\}=\sum _{l=1}^{l}C_{i,j,k}^lv_l$$
and satisfy
$$C_{ijk}^l=(-1)^{\varepsilon(\sigma)}C_{\sigma(i)\sigma(j)\sigma(k)}^l,$$
for any $\sigma \in \Sigma_3.$ The Jacobi condition writes
$$
\begin{array}{l}
C_{ijk}^t C_{tlm}^s-C_{ijl}^t C_{tkm}^s+C_{ijm}^tC_{tjk}^s+C_{ikl}^tC_{tjm}^s-C_{ikm}^tC_{tjl}^s+C_{ilm}^tC_{tjl}^s\\
-C_{jkl}^tC_{tim}^s+C_{jkm}^tC_{til}^s-C_{jlm}^tC_{tik}^s+C_{klm}^tC_{tij}^s=0,
\end{array}
$$
for any $i<j<k,l<m$ and $s,t=1,\cdots,p.$ Let $\{\omega_1,\cdots,\omega_p\}$ be the dual basis of $\{v_1,\cdots,v_p\}$. We consider the graded exterior algebra $\Lambda (V)=\oplus \Lambda ^{k}$ of $V$ and the linear operator
$$d: \Lambda ^1(V)=V^*\rightarrow \Lambda ^3 (V)$$
given by
$$d\omega_l=\sum_{i<j<k} C_{ijk}^l\omega_i\wedge\omega_j\wedge\omega_k.$$
If we denote also by $d$ the linear operator
$$d: \Lambda ^3(V)=V^*\rightarrow \Lambda ^5 (V)$$
defined by
$$d(\omega_i\wedge\omega_j\wedge\omega_k)=d\omega_i\wedge\omega_j\wedge\omega_k+\omega_i\wedge d\omega_j\wedge\omega_k+\omega_i\wedge\omega_j\wedge d\omega_k,$$
we obtain
$$
\begin{array}{ll}
d(d\omega_l)&=\displaystyle\sum_{i<j<k} C_{ijk}^l ( d\omega_i\wedge\omega_j\wedge\omega_k+\omega_i\wedge d\omega_j\wedge\omega_k+\omega_i\wedge\omega_j\wedge d\omega_k \\
&= \displaystyle\sum_{i<j<k}C_{ijk}^l(C_{lst}^i\omega_l\wedge\omega_s\wedge\omega_t\wedge\omega_j\wedge\omega_k +C_{lst}^j \omega_i\wedge\omega_l\wedge\omega_s\wedge\omega_t\wedge\omega_k\\
& \ \  \ \  +C_{lst}^k \omega_i\wedge\omega_j\wedge\omega_l\wedge\omega_s\wedge\omega_t).
\end{array}
$$
In this summand, all the products containing two equal factors  are zero (this justifies to use  the exterior algebra). In the same way, the Jacobi condition related to five vectors is trivial as soon as two vectors are equal. In fact, if we compute the Jacobi condition for the vectors $(X_1,X_2,X_3,X_1,X_1)$ we find $0=0$ and for the vector $(X_1,X_2,X_3,X_1,X_5)$ we find
$$
\begin{array}{lll}
[[X_1,X_2,X_3],X_1,X_5]+[[X_1,X_2,X_5],X_3,X_1]-[[X_1,X_3,X_5],X_2,X_1]&&\\
-[[X_2,X_3,X_1],X_1,X_5]-[[X_2,X_1,X_5],X_1,X_3]-[[X_3,X_1,X_5],X_1,X_2]&=&0,
\end{array}$$
that is, $0=0$. Thus the Jacobi condition concerns a  family of $5$ independent vectors. Let us return to the computation of $d(d\omega)$. The coefficient of $d(d\omega _l)$ related for example to
$\omega_1\wedge\omega_2\wedge\omega_3\wedge\omega_4\wedge\omega_5$ corresponds to the coefficient of $X_l$ in the Jacobi condition related to $(X_1,X_2,X_3,X_4,X_5)$. Thus
$$d(d\omega_l)=0.$$
These relations can be called the Maurer-Cartan equations.

\medskip

\noindent {\bf Remark.} We cannot use the same calculus to obtain Maurer-Cartan equations adapted to the structure of $n$-Lie algebras. This means that the Maurer-Cartan equations of a $n$-Lie algebra are the Maurer-Cartan equations of this algebra considered as a sh-$n$-Lie algebra. In the classical case of Lie algebras, we have also  such a situation. For example, when we consider the $2$-step nilpotent Lie algebras, defined by the $2$-step Jacobi condition
$$[[X_i,X_j],X_k]=0,$$
there is no  exterior calculus adapted to this special Jacobi condition. The Maurer-Cartan equations of a $2$-step nilpotent algebra are the Maurer-Cartan equations of this algebra considered as a Lie algebra.

\end{document}